\documentclass{amsart}

\usepackage{xypic}
\input xy
\xyoption{all}
\usepackage{epsfig}
\usepackage{amsthm}
\usepackage{amssymb}
\usepackage{amsmath}
\usepackage{amscd}
\usepackage{color}
\usepackage[T1]{fontenc}

\newcommand{\bg}{\begin{equation}}
\newcommand{\ed}{\end{equation}}
\newcommand{\bga}{\begin{eqnarray}}
\newcommand{\eda}{\end{eqnarray}}
\newcommand{\pf}{\textbf{Proof:\ }}

\def\cbdu{\par{\raggedleft$\Box$\par}}

\newtheorem {Theorem}  {Theorem}

\numberwithin{Theorem}{section}

\newtheorem {Lemma}[Theorem]  {Lemma}

\theoremstyle{definition}

\theoremstyle{remark}
\newtheorem{Remark}[Theorem]{\bf Remark}

\expandafter\chardef\csname pre amssym.def
at\endcsname=\the\catcode`\@ \catcode`\@=11
\def\undefine#1{\let#1\undefined}
\def\newsymbol#1#2#3#4#5{\let\next@\relax
 \ifnum#2=\@ne\let\next@\msafam@\else
 \ifnum#2=\tw@\let\next@\msbfam@\fi\fi
 \mathchardef#1="#3\next@#4#5}
\def\mathhexbox@#1#2#3{\relax
 \ifmmode\mathpalette{}{\m@th\mathchar"#1#2#3}%
 \else\leavevmode\hbox{$\m@th\mathchar"#1#2#3$}\fi}
\def\hexnumber@#1{\ifcase#1 0\or 1\or 2\or 3\or 4\or 5\or 6\or 7\or 8\or
 9\or A\or B\or C\or D\or E\or F\fi}

\font\teneufm=eufm10 \font\seveneufm=eufm7 \font\fiveeufm=eufm5
\newfam\eufmfam
\textfont\eufmfam=\teneufm \scriptfont\eufmfam=\seveneufm
\scriptscriptfont\eufmfam=\fiveeufm

\catcode`\@=\csname pre amssym.def at\endcsname

\newcounter{remark}
\newcommand{\R}{\mathbf{R}}

\def \ls {{\lambda_q^{2s}}}

\def  \R   {{\mathbb R}}

\def  \12  {{\frac{1}{2}}}

\def\build#1_#2^#3{\mathrel{\mathop{\kern 0pt#1}\limits_{#2}^{#3}}}

\begin{document}

\title[Local well-posedness of Hall-MHD]{Local well-posedness of the Hall-MHD system in $H^s(\R^n)$ with $s>\frac n2$}


\author [Mimi Dai]{Mimi Dai}
\address{Department of Mathematics, Stat. and Comp.Sci., University of Illinois Chicago, Chicago, IL 60607,USA}
\email{mdai@uic.edu}





\begin{abstract}
We establish local well-posedness of the Hall-magneto-hydrodynamics (Hall-MHD) system in the Sobolev space $\left(H^s(\R^n)\right)^2$ with $s>\frac n2$. The previously known local well-posedness space was $\left(H^s(\R^n)\right)^2$ with $s>\frac n2+1$. Thus the result presented here is an improvement.

\bigskip

KEY WORDS: Hall-magneto-hydrodynamics; local well-posedness; Littlewood-Paley theory.

\hspace{0.02cm}CLASSIFICATION CODE: 76D03, 76W05, 35Q35, 35D35.
\end{abstract}

\maketitle

\section{Introduction}

Considered here is the incompressible 
Hall-magneto-hydrodynamics (Hall-MHD) system with fractional magnetic diffusion:
\begin{equation}\label{HMHD}
\begin{split}
u_t+u\cdot\nabla u-b\cdot\nabla b+\nabla p=\nu\Delta u,\\
b_t+u\cdot\nabla b-b\cdot\nabla u+\eta\nabla\times((\nabla\times b)\times b)=-\mu(-\Delta)^\alpha b,\\
\nabla \cdot u=0, 
\end{split}
\end{equation}
with  $(x,t)\in\mathbb{R}^n\times [0,\infty)$, $n\geq 2$, 
and initial conditions
\begin{equation}
u(x,0)=u_0(x),\qquad b(x,0)=b_0(x), \qquad \nabla\cdot u_0=\nabla\cdot b_0=0.
\end{equation}
Here $u$ is the fluid velocity, $p$ pressure and $b$ the magnetic field.  The constants $\nu,\mu$ and $\eta$ denote the kinematic viscosity, the reciprocal of the magnetic Reynolds number and the Hall effect coefficient, respectively. We assume $\nu>0$, $\mu>0$ and $\alpha>\frac12$.
The Hall term $\nabla\times((\nabla\times b)\times b)$ is the only difference between the Hall-MHD  and the usual MHD system. For
mathematical study on this model, we refer to  \cite{ADFL, CL, CS, CWW, CW, D, DL, DS} and reference therein.


The purpose of this paper is to find the largest possible Sobolev spaces where the Hall-MHD system is locally well-posed. Previously, it was shown in \cite{CW} that system (\ref{HMHD}) with $\alpha=1$ is locally well-posed in $\left(H^s(\R^3)\right)^2$ with $s>\frac 52$. Later, in the case of $\frac12<\alpha<1$, local well-posedness was obtained in $\left(H^s(\R^n)\right)^2$ with $s>\frac n2+1$. We aim to improve the aforementioned findings and establish the main result below.

\begin{Theorem}\label{thm}
Let $\nu, \mu> 0$ and $\alpha> \frac12$. Assume $(u_0,b_0)\in \left(H^s(\R^n)\right)^2$ with $s>2-2\alpha+\frac{n}2$ and $\nabla \cdot u_0=\nabla\cdot b_0=0$. There exists a time $T=T(\|u_0\|_{H^s},\|b_0\|_{H^s})>0$ and a unique solution $(u,b)$ of (\ref{HMHD}) on $[0,T]$ such that 
\[(u,b)\in \left(C([0,T]; H^s(\R^n))\right)^2.\]
\end{Theorem}


\begin{Remark}\label{remk}
 Notice that $s>2-2\alpha+\frac{n}2=\frac n2$ for $\alpha=1$; and $2-2\alpha+\frac{n}2<\frac n2+1$ for $\frac 12<\alpha<1$. Thus for the regular Hall-MHD system, that is (\ref{HMHD}) with $\alpha=1$, we obtain the local well-posedness in $\left(H^s(\R^n)\right)^2$ with $s>\frac n2$, which is a lager than $\left(H^{5/2}(\R^3)\right)^2$ for $n=3$.
\end{Remark}

The techniques involved are based on the Littlewood-Paley decomposition theory and the frequency-localization approach.

\medskip


\noindent
\textbf{Notation.}
For the sake of brevity, we denote by:
 $A\lesssim B$ an estimate of the form $A\leq C B$ with
an absolute constant $C$; $A\sim B$ an estimate of the form $C_1
B\leq A\leq C_2 B$ with absolute constants $C_1$, $C_2$;
 $\|\cdot\|_p$ the norm of space $L^p$; 
 and $(\cdot, \cdot)$ the $L^2$-inner product. The notations associated with Littlewood-Paley decomposition theory and related concepts are introduced in Appendix.

\section{A priori estimate}
\label{sec-est}

The core of the proof of local well-posedness is the a priori estimate satisfied by smooth solutions in $H^s$ with $s>2-2\alpha+\frac n2$, which is the content of this section. The local existence of smooth solutions will then follow from certain traditional approximating and limiting process. The uniqueness and continuous dependance on initial data can be also obtained through standard arguments. 
Thus, we only show

\begin{Theorem}\label{le-priori}
Let $(u_0, b_0)\in (H^s(\R^n))^2$ with $s>2-2\alpha+\frac n2$ and $(u, b)$ be a smooth solution of (\ref{HMHD}) starting from the data $(u_0, b_0)$. There exists a time $T=T(\|u_0\|_{H^s}, \|b_0\|_{H^s})> 0$, such that, for every $t\in [0,T]$ we have
\begin{equation}\notag
\|u(t)\|_{H^s}^2+\|b(t)\|_{H^s}^2 \leq C\left( \|u_0\|_{H^s}^2+\|b_0\|_{H^s}^2\right),
\end{equation}
where the constant $C$ depends on $T$, $\nu$, $\mu$, $\|u_0\|_{H^s}$, and $\|b_0\|_{H^s}$.
\end{Theorem}
\pf
Multiplying the first equation of (\ref{HMHD}) by $\lambda_q^{2s}\Delta_q^2 u$ and the second one by $\lambda_q^{2s}\Delta_q^2 b$, and taking summation for all $q\geq -1$ gives us  
\begin{equation}\label{ineq-ubq}
\begin{split}
&\frac12\frac{d}{dt}\sum_{q\geq -1}\left(\lambda_q^{2s}\|u_q\|_2^2+\lambda_q^{2s}\|b_q\|_2^2\right)\\
\leq &-\nu\sum_{q\geq-1}\lambda_q^{2s+2}\|u_q\|_2^2-\mu\sum_{q\geq-1}\lambda_q^{2s+2\alpha}\|b_q\|_2^2+ I_1+I_2+I_3+I_4+I_5,
\end{split}
\end{equation}
with
\begin{equation}\notag
\begin{split}
I_1=&-\sum_{q\geq -1}\lambda_q^{2s}\int_{\R^3}\Delta_q(u\cdot\nabla u)\cdot u_q\, dx, \qquad
I_2=\sum_{q\geq -1}\lambda_q^{2s}\int_{\R^3}\Delta_q(b\cdot\nabla b)\cdot u_q\, dx,\\
I_3=&-\sum_{q\geq -1}\lambda_q^{2s}\int_{\R^3}\Delta_q(u\cdot\nabla b)\cdot b_q\, dx,\qquad
I_4=\sum_{q\geq -1}\lambda_q^{2s}\int_{\R^3}\Delta_q(b\cdot\nabla u)\cdot b_q\, dx,\\
I_5=&\sum_{q\geq -1}\lambda_q^{2s}\int_{\R^3}\Delta_q((\nabla\times b)\times b)\cdot \nabla\times b_q\, dx.
\end{split}
\end{equation}
As expected, the estimate of $I_1, I_2, I_3$, and $I_4$ are less challenging than that of $I_5$. On the other hand, due to the similarity of $I_1$ and $I_3$, $I_2$ and $I_4$, we are eligible to only show the details of handling $I_3$ and $I_2$, not $I_1$ and $I_4$. 

We first decompose $I_3$ by adapting Bony's paraproduct (\ref{Bony})
\begin{equation}\notag
\begin{split}
I_3=
&-\sum_{q\geq -1}\sum_{|q-p|\leq 2}\lambda_q^{2s}\int_{\R^3}\Delta_q(u_{\leq p-2}\cdot\nabla b_p)\cdot b_q\, dx\\
&-\sum_{q\geq -1}\sum_{|q-p|\leq 2}\lambda_q^{2s}\int_{\R^3}\Delta_q(u_{p}\cdot\nabla b_{\leq{p-2}})\cdot b_q\, dx\\
&-\sum_{q\geq -1}\sum_{p\geq q-2}\lambda_q^{2s}\int_{\R^3}\Delta_q(u_p\cdot\nabla\tilde b_p)\cdot b_q\, dx\\
=&I_{31}+I_{32}+I_{33};
\end{split}
\end{equation}
and then by commutator (\ref{commu}) to rewrite $I_{31}$
\begin{equation}\notag
\begin{split}
I_{31}=&-\sum_{q\geq -1}\sum_{|q-p|\leq 2}\lambda_q^{2s}\int_{\R^3}[\Delta_q, u_{\leq{p-2}}\cdot\nabla] b_p\cdot b_q\, dx\\
&-\sum_{q\geq -1}\sum_{|q-p|\leq 2}\lambda_q^{2s}\int_{\R^3}(u_{\leq{q-2}}\cdot\nabla \Delta_q b_p)\cdot  b_q\, dx\\
&-\sum_{q\geq -1}\sum_{|q-p|\leq 2}\lambda_q^{2s}\int_{\R^3}\left((u_{\leq{p-2}}-u_{\leq{q-2}})\cdot\nabla\Delta_qb_p\right)\cdot b_q\, dx\\
=&I_{311}+I_{312}+I_{313}.
\end{split}
\end{equation}
Since $\sum_{|p-q|\leq 2}\Delta_qb_p=b_q$ and $\nabla\cdot u_{\leq q-2}=0$, one can infer $I_{312}=0$.  

To estimate $I_{311}$, it follows from the commutator estimate in Lemma \ref{le-commu}, H\"older's inequality, and Bernstein's inequality that
\begin{equation}\notag
\begin{split}
|I_{311}|\leq&\sum_{q\geq -1}\sum_{|p-q|\leq 2}\lambda_q^{2s}\|\nabla u_{\leq p-2}\|_\infty\|b_p\|_2\|b_q\|_2\\
\lesssim &\sum_{q\geq -1}\lambda_q^{2s}\|b_q\|_2^2\sum_{p\leq q}\lambda_p^{1+\frac n2}\|u_p\|_2\\
\lesssim&\sum_{q\geq -1}\sum_{p\leq q}\lambda_{p-q}^{\delta\alpha}\lambda_p^{1+\frac n2-\delta-\delta\alpha-s}\left(\lambda_q^{s+\alpha}\|b_q\|_2\right)^\delta\left(\lambda_q^{s}\|b_q\|_2\right)^{2-\delta}\\
&\cdot\left(\lambda_p^{s+1}\|u_p\|_2\right)^\delta\left(\lambda_p^{s}\|u_p\|_2\right)^{1-\delta}\\
\lesssim & \sum_{q\geq -1}\sum_{p\leq q}\lambda_{p-q}^{\delta\alpha}\left(\lambda_q^{s+\alpha}\|b_q\|_2\right)^\delta\left(\lambda_q^{s}\|b_q\|_2\right)^{2-\delta}\left(\lambda_p^{s+1}\|u_p\|_2\right)^\delta\left(\lambda_p^{s}\|u_p\|_2\right)^{1-\delta}
\end{split}
\end{equation}
for some parameter $0<\delta<1$ satisfying 
\begin{equation}\label{para3}
s\geq 1+\frac n2-\delta-\delta\alpha.
\end{equation}
We continue the estimate of $I_{311}$ by using Young's inequality with parameters satisfying 
\begin{equation}\label{para4}
\begin{split}
\frac1{\delta_1}+\frac1{\delta_2}+\frac1{\delta_3}+\frac1{\delta_4}=\delta\alpha, \ \ 0<\delta_1,\delta_2,\delta_3,\delta_4<1\\
\frac1{\theta_1}+\frac1{\theta_2}+\frac1{\theta_3}+\frac1{\theta_4}=1, \ \ \theta_1=\theta_3=\frac2{\delta}, \ \ 1<\theta_2, \theta_4<\infty.
\end{split}
\end{equation}
It then follows that 
\begin{equation}\notag
\begin{split}
|I_{311}|
\leq & \frac{\mu}{16}\sum_{q\geq -1}\sum_{p\leq q}\lambda_{p-q}^{\delta_1\theta_1}\lambda_q^{2s+2\alpha}\|b_q\|_2^2+C_{\nu,\mu}\sum_{q\geq -1}\sum_{p\leq q}\lambda_{p-q}^{\delta_2\theta_2}\left(\lambda_q^{s}\|b_q\|_2\right)^{(2-\delta)\theta_2}\\
&+\frac\nu{16}\sum_{q\geq -1}\sum_{p\leq q}\lambda_{p-q}^{\delta_3\theta_3}\lambda_p^{2s+2}\|u_p\|_2^2+C_{\nu,\mu}\sum_{q\geq -1}\sum_{p\leq q}\lambda_{p-q}^{\delta_4\theta_4}\left(\lambda_q^{s}\|u_p\|_2\right)^{(1-\delta)\theta_4}\\
\leq & \frac\mu{16}\sum_{q\geq -1}\lambda_q^{2s+2\alpha}\|b_q\|_2^2+ \frac\nu{16}\sum_{q\geq -1}\lambda_q^{2s+2}\|u_q\|_2^2\\
&+C_{\nu,\mu}\sum_{q\geq -1}\left(\lambda_q^{s}\|b_q\|_2\right)^{(2-\delta)\theta_2}+C_{\nu,\mu}\sum_{q\geq -1}\left(\lambda_q^{s}\|u_q\|_2\right)^{(1-\delta)\theta_4},
\end{split}
\end{equation}
with various constants $C_{\nu,\mu}$ that depend on $\nu,\mu$ and tend to infinity as $\nu,\mu\to 0$.
We pause to analyze the parameters. In view of (\ref{para3}) and (\ref{para4}), we obtain that
\begin{equation}\label{s-i3}
s\geq \frac n2-\alpha+\left(\frac1{\theta_2}+\frac1{\theta_4}\right)(1+\alpha)\geq\frac n2-\alpha+\epsilon
\end{equation}
provided $\theta_2$ and $\theta_4$ are large enough. 

Other terms in $I_3$ are simpler and can be estimated in an analogous way; thus the details are omitted.
As a conclusion, we have for $s$ satisfying (\ref{s-i3})
\begin{equation}\label{est-i3}
\begin{split}
|I_3|\leq &\frac\mu{8}\sum_{q\geq -1}\lambda_q^{2s+2\alpha}\|b_q\|_2^2+ \frac\nu{8}\sum_{q\geq -1}\lambda_q^{2s+2}\|u_q\|_2^2\\
&+C_{\nu,\mu}\left(\sum_{q\geq -1}\lambda_q^{2s}\|b_q\|_2^2\right)^{\gamma_1}+C_{\nu,\mu}\left(\sum_{q\geq -1}\lambda_q^{2s}\|u_q\|_2^2\right)^{\gamma_2},
\end{split}
\end{equation}
with certain constants $\gamma_1,\gamma_2>1$.

Adapting the same decomposition strategy of using Bony's paraproduct and commutator,  we deconstruct $I_2$ and $I_4$ as follows
\begin{equation}\notag
\begin{split}
I_2=
&\sum_{q\geq -1}\sum_{|q-p|\leq 2}\lambda_q^{2s}\int_{\R^3}\Delta_q(b_{\leq p-2}\cdot\nabla b_p)\cdot u_q\, dx\\
&+\sum_{q\geq -1}\sum_{|q-p|\leq 2}\lambda_q^{2s}\int_{\R^3}\Delta_q(b_{p}\cdot\nabla b_{\leq{p-2}})\cdot u_q\, dx\\
&+\sum_{q\geq -1}\sum_{p\geq q-2}\lambda_q^{2s}\int_{\R^3}\Delta_q(b_p\cdot\nabla\tilde b_p)\cdot u_q\, dx\\
=&I_{21}+I_{22}+I_{23},
\end{split}
\end{equation}
with 
\begin{equation}\notag
\begin{split}
I_{21}=&\sum_{q\geq -1}\sum_{|q-p|\leq 2}\lambda_q^{2s}\int_{\R^3}[\Delta_q, b_{\leq{p-2}}\cdot\nabla] b_p\cdot u_q\, dx\\
&+\sum_{q\geq -1}\sum_{|q-p|\leq 2}\lambda_q^{2s}\int_{\R^3}(b_{\leq{q-2}}\cdot\nabla \Delta_q b_p)\cdot u_q\, dx\\
&+\sum_{q\geq -1}\sum_{|q-p|\leq 2}\lambda_q^{2s}\int_{\R^3}((b_{\leq{p-2}}-b_{\leq{q-2}})\cdot\nabla\Delta_qb_p)\cdot u_q\, dx\\
=&I_{211}+I_{212}+I_{213};
\end{split}
\end{equation}
and
\begin{equation}\notag
\begin{split}
I_4=
&\sum_{q\geq -1}\sum_{|q-p|\leq 2}\lambda_q^{2s}\int_{\R^3}\Delta_q(b_{\leq p-2}\cdot \nabla u_p) \cdot b_q\, dx\\
&+\sum_{q\geq -1}\sum_{|q-p|\leq 2}\lambda_q^{2s}\int_{\R^3}\Delta_q(b_{p}\cdot \nabla u_{\leq{p-2}}) \cdot b_q\, dx\\
&+\sum_{q\geq -1}\sum_{p\geq q-2}\lambda_q^{2s}\int_{\R^3}\Delta_q(\tilde b_p\cdot \nabla u_p)\cdot b_q\, dx\\
=&I_{41}+I_{42}+I_{43},
\end{split}
\end{equation}
with 
\begin{equation}\notag
\begin{split}
I_{41}=&\sum_{q\geq -1}\sum_{|q-p|\leq 2}\lambda_q^{2s}\int_{\R^3}[\Delta_q, b_{\leq{p-2}}\cdot\nabla] u_p\cdot b_q\, dx\\
&+\sum_{q\geq -1}\sum_{|q-p|\leq 2}\lambda_q^{2s}\int_{\R^3}(b_{\leq{q-2}}\cdot\nabla \Delta_q u_p)\cdot b_q\, dx\\
&+\sum_{q\geq -1}\sum_{|q-p|\leq 2}\lambda_q^{2s}\int_{\R^3}((b_{\leq{p-2}}-b_{\leq{q-2}})\cdot\nabla\Delta_qu_p)\cdot b_q\, dx\\
=&I_{411}+I_{412}+I_{413}.
\end{split}
\end{equation}
We claim that $I_{212}+I_{412}=0$. Indeed, we have
\begin{equation}\notag
\begin{split}
I_{212}+I_{412}
=&\sum_{q\geq -1}\sum_{|q-p|\leq 2}\lambda_q^{2s}\int_{\R^3}(b_{\leq{q-2}}\cdot\nabla \Delta_q b_p)\cdot u_q\, dx\\
&+\sum_{q\geq -1}\sum_{|q-p|\leq 2}\lambda_q^{2s}\int_{\R^3}(b_{\leq{q-2}}\cdot\nabla \Delta_q u_p)\cdot b_q\, dx\\
=&\sum_{q\geq -1}\lambda_q^{2s}\int_{\R^3}(b_{\leq{q-2}}\cdot\nabla  b_q)\cdot u_q\, dx
+\sum_{q\geq -1}\lambda_q^{2s}\int_{\R^3}(b_{\leq{q-2}}\cdot\nabla u_q)\cdot b_q\, dx\\
=&0.
\end{split}
\end{equation}
The fact $\sum_{|p-q|\leq 2}\Delta_qb_p=b_q$ and $\sum_{|p-q|\leq 2}\Delta_qu_p=u_q$ justifies the second equality above. 

The rest terms in $I_2+I_4$ are relatively simple. We only choose one representative term, $I_{211}$, to carry out the details of estimating.  Applying H\"older's inequality and Bernstein's inequality leads to
\begin{equation}\notag
\begin{split}
|I_{211}|\leq & \sum_{q\geq -1}\sum_{|q-p|\leq 2}\lambda_q^{2s}\|\nabla b_{\leq p-2}\|_\infty\|b_p\|_2\|u_q\|_2\\
\lesssim & \sum_{q\geq -1}\lambda_q^{2s}\|b_q\|_2\|u_q\|_2\sum_{p\leq q}\lambda_p^{1+\frac n2}\|b_p\|_2\\
=&\sum_{q\geq -1}\sum_{p\leq q}\lambda_{p-q}^{\delta_1\alpha+\delta_2}\lambda_p^{1+\frac n2-\delta_1\alpha-\delta_3\alpha-\delta_2-s}\left(\lambda_q^{s+\alpha}\|b_q\|_2\right)^{\delta_1}
\left(\lambda_q^{s}\|b_q\|_2\right)^{1-\delta_1}\\
&\cdot \left(\lambda_q^{s+1}\|u_q\|_2\right)^{\delta_2}
\left(\lambda_q^{s}\|u_q\|_2\right)^{1-\delta_2}\left(\lambda_p^{s+\alpha}\|b_p\|_2\right)^{\delta_3}
\left(\lambda_p^{s}\|b_p\|_2\right)^{1-\delta_3}\\
\leq &C\sum_{q\geq -1}\sum_{p\leq q}\lambda_{p-q}^{\delta_1\alpha+\delta_2}\left(\lambda_q^{s+\alpha}\|b_q\|_2\right)^{\delta_1}
\left(\lambda_q^{s}\|b_q\|_2\right)^{1-\delta_1}
\cdot \left(\lambda_q^{s+1}\|u_q\|_2\right)^{\delta_2}\\
&\left(\lambda_q^{s}\|u_q\|_2\right)^{1-\delta_2}
\left(\lambda_p^{s+\alpha}\|b_p\|_2\right)^{\delta_3}
\left(\lambda_p^{s}\|b_p\|_2\right)^{1-\delta_3}\\
\end{split}
\end{equation}
for parameters $0<\delta_1, \delta_2,\delta_3<1$, $\delta_2=(2-\delta_1-\delta_2)\alpha$, and 
\begin{equation}\label{para241}
s\geq 1+\frac n2-\delta_1\alpha-\delta_3\alpha-\delta_2.
\end{equation} 
Adapting Young's inequality with parameters $\zeta_i$, $1\leq i\leq 6$, such that 
\begin{equation}\label{para242}
\begin{split}
\zeta_1+\zeta_2+\zeta_3+\zeta_4+\zeta_5+\zeta_6=\delta_1\alpha+\delta_2,  \ \ \zeta_1, ..., \zeta_6>0\\
\frac1{\theta_1}+\frac1{\theta_2}+\frac1{\theta_3}+\frac1{\theta_4}+\frac1{\theta_5}+\frac1{\theta_6}=1, \\
 \theta_1=\frac2{\delta_1}, \ \theta_3=\frac2{\delta_2}, \ \theta_5=\frac2{\delta_3}, \ 1<\theta_2,\theta_4,\theta_6<\infty
\end{split}
\end{equation}
we have
\begin{equation}\notag
\begin{split}
|I_{211}|\leq &\frac{\mu}{8}\sum_{g\geq -1}\lambda_q^{2s+2\alpha}\|b_q\|_2^2+\frac{\nu}{16}\sum_{g\geq -1}\lambda_q^{2s+2}\|u_q\|_2^2+C_{\nu,\mu}\sum_{g\geq -1}\left(\lambda_q^{s}\|b_q\|_2\right)^{(1-\delta_1)\theta_2}\\
&+C_{\nu,\mu}\sum_{g\geq -1}\left(\lambda_q^{s}\|b_q\|_2\right)^{(1-\delta_3)\theta_6}
+C_{\nu,\mu}\sum_{g\geq -1}\left(\lambda_q^{s}\|u_q\|_2\right)^{(1-\delta_2)\theta_4}.
\end{split}
\end{equation}
Again, the parameter constraints (\ref{para241}) and (\ref{para242}) imply that
\begin{equation}\notag
\begin{split}
s\geq &1+\frac n2-2\alpha+(\alpha-1)\delta_2+2\alpha\left(\frac1{\theta_2}+\frac1{\theta_4}+\frac1{\theta_6}\right)\\
=&1+\frac n2-2\alpha+(\alpha-1)\delta_2+\epsilon
\end{split}
\end{equation}
for large enough $\theta_2,\theta_4$, and $\theta_6$. Notice that $s\geq \frac n2-1+\epsilon$ for $\alpha=1$. In general for $\delta_2$ close enough to 1, we have
\begin{equation}\label{s-i2}
s\geq \frac n2-\alpha+\epsilon.
\end{equation}
To conclude, we expect to have for $s$ satisfying (\ref{s-i2})
\begin{equation}\label{est-i2}
\begin{split}
|I_{2}|\leq &\frac{\mu}{8}\sum_{g\geq -1}\lambda_q^{2s+2\alpha}\|b_q\|_2^2+\frac{\nu}{16}\sum_{g\geq -1}\lambda_q^{2s+2}\|u_q\|_2^2\\
&+C_{\nu,\mu}\left(\sum_{q\geq -1}\lambda_q^{2s}\|b_q\|_2^2\right)^{\gamma_1}+C_{\nu,\mu}\left(\sum_{q\geq -1}\lambda_q^{2s}\|u_q\|_2^2\right)^{\gamma_2},
\end{split}
\end{equation}
for some constants $\gamma_1,\gamma_2$.

Now we are left to estimate $I_5$. By Bony's paraproduct and commutator (\ref{comm-v}), the routine decomposition procedure yields
\begin{equation}\notag
\begin{split}
I_{5}=&\sum_{q\geq-1}\sum_{|q-p|\leq 2}\lambda_q^{2s}\int_{\mathbb R^3}\Delta_q( b_{\leq p-2}\times(\nabla\times b_p))\cdot\nabla\times b_q\, dx\\
&+\sum_{q\geq-1}\sum_{|q-p|\leq 2}\lambda_q^{2s}\int_{\mathbb R^3}\Delta_q( b_{p}\times(\nabla\times b_{\leq p-2}))\cdot\nabla\times b_q\, dx\\
&+\sum_{q\geq-1}\sum_{p\geq q-2}\lambda_q^{2s}\int_{\mathbb R^3}\Delta_q( b_{p}\times(\nabla\times \tilde b_p))\cdot\nabla\times b_q\, dx\\
=&I_{51}+I_{52}+I_{53};
\end{split}
\end{equation}
with 
\begin{equation}\notag
\begin{split}
I_{51}=&\sum_{q\geq-1}\sum_{|q-p|\leq 2}\lambda_q^{2s}\int_{\mathbb R^3}[\Delta_q,b_{\leq p-2}\times\nabla\times]b_p\cdot\nabla\times b_q\, dx\\
&+\sum_{q\geq-1}\lambda_q^{2s}\int_{\mathbb R^3}b_{\leq q-2}\times(\nabla\times b_q)\cdot\nabla\times b_q\, dx\\
&+\sum_{q\geq-1}\sum_{|p-q|\leq 2}\lambda_q^{2s}\int_{\mathbb R^3}(b_{\leq p-2}-b_{\leq q-2})\times(\nabla\times (b_p)_q)\cdot\nabla\times b_q\, dx\\
=&I_{511}+I_{512}+I_{513}.
\end{split}
\end{equation}
The cross product property implies immediately that $I_{512}=0$. 
We deduce from the commutator estimate in Lemma  \ref{le-Hall1} that
\begin{equation}\notag
\begin{split}
|I_{511}|\lesssim &\sum_{q\geq-1}\sum_{|p-q|\leq 2}\lambda_q^{2s+1}\|\nabla b_{\leq p-2}\|_\infty\|b_p\|_2\|b_q\|_2\\
\lesssim & \sum_{q\geq -1}\lambda_q^{2s+1}\|b_q\|_2^2\sum_{p\leq q}\lambda_p^{1+\frac n2}\|b_p\|_2\\
=&\sum_{q\geq -1}\sum_{p\leq q}\lambda_{p-q}^{\delta_1\alpha-1}\lambda_p^{2+\frac n2-\delta_1\alpha-\delta_2\alpha-s}\left(\lambda_q^{s+\alpha}\|b_q\|_2\right)^{\delta_1}\left(\lambda_q^{s}\|b_q\|_2\right)^{2-\delta_1}\\
&\cdot\left(\lambda_p^{s+\alpha}\|b_p\|_2\right)^{\delta_2}\left(\lambda_p^{s}\|b_p\|_2\right)^{1-\delta_2}\\
\leq &C \sum_{q\geq -1}\sum_{p\leq q}\lambda_{p-q}^{\delta_1\alpha-1}\left(\lambda_q^{s+\alpha}\|b_q\|_2\right)^{\delta_1}\left(\lambda_q^{s}\|b_q\|_2\right)^{2-\delta_1}\left(\lambda_p^{s+\alpha}\|b_p\|_2\right)^{\delta_2}\left(\lambda_p^{s}\|b_p\|_2\right)^{1-\delta_2}
\end{split}
\end{equation}
for parameters satisfying $\frac 1\alpha<\delta_1<2$, $0<\delta_2<1$, and 
\begin{equation}\label{para5}
s\geq 2+\frac n2-\delta_1\alpha-\delta_2\alpha.
\end{equation}
By Young's inequality we have for the parameters 
\begin{equation}\label{para6}
\begin{split}
\zeta_1+\zeta_2+\zeta_3+\zeta_4=\delta_1\alpha-1, \ \ \zeta_1, ..., \zeta_4>0\\
\frac1{\theta_1}+\frac1{\theta_2}+\frac1{\theta_3}+\frac1{\theta_4}=1, \ \ \theta_1=\frac2{\delta_1}, \ \ \theta_3=\frac2{\delta_2}, \ \ 1<\theta_2, \theta_4<\infty,
\end{split}
\end{equation}
such that
\begin{equation}\notag
\begin{split}
|I_{511}|\leq & \frac{\mu}{16}\sum_{q\geq -1}\sum_{p\leq q}\lambda_{p-q}^{\zeta_1\theta_1}\lambda_q^{2s+2\alpha}\|b_q\|_2^2+C_\mu\sum_{q\geq -1}\sum_{p\leq q}\lambda_{p-q}^{\zeta_2\theta_2}\left(\lambda_q^{s}\|b_q\|_2\right)^{(2-\delta_1)\theta_2}\\
&+\frac{\mu}{16}\sum_{q\geq -1}\sum_{p\leq q}\lambda_{p-q}^{\zeta_3\theta_3}\lambda_p^{2s+2\alpha}\|b_p\|_2^2+C_\mu\sum_{q\geq -1}\sum_{p\leq q}\lambda_{p-q}^{\zeta_4\theta_4}\left(\lambda_p^{s}\|b_p\|_2\right)^{(1-\delta_2)\theta_4}\\
\leq & \frac{\mu}{8}\sum_{q\geq -1}\lambda_q^{2s+2\alpha}\|b_q\|_2^2+C_\mu\sum_{q\geq -1}\left(\lambda_q^{s}\|b_q\|_2\right)^{(2-\delta_1)\theta_2}\\
&+C_\mu\sum_{q\geq -1}\left(\lambda_p^{s}\|b_p\|_2\right)^{(1-\delta_2)\theta_4}.
\end{split}
\end{equation}
Regarding the parameters, (\ref{para5}) and (\ref{para6}) imply that
\begin{equation}\label{s-i51}
s\geq \frac n2+2-2\alpha+2\alpha\left(\frac1{\theta_2}+\frac1{\theta_4}\right)\geq\frac n2+2-2\alpha+\epsilon
\end{equation}
for large enough $\theta_2$ and $\theta_4$.

By H\"older's inequality,
\begin{equation}\notag
\begin{split}
|I_{513}|\leq &\sum_{q\geq-1}\sum_{|p-q|\leq 2}\lambda_q^{2s}\int_{\mathbb R^3}\left|(b_{\leq p-2}-b_{\leq q-2})\times(\nabla\times (b_p)_q)\cdot\nabla\times b_q\right|\, dx\\
\lesssim &\sum_{q\geq -1}\sum_{|p-q|\leq 2}\lambda_q^{2s}\|\nabla b_q\|_\infty\|b_{\leq p-2}-b_{\leq q-2}\|_2\|\nabla b_p\|_2\\
\lesssim &\sum_{q\geq -1}\lambda_q^{2s+2+\frac n2}\|b_q\|_2^3\\
=& C\sum_{q\geq -1}\lambda_q^{2+\frac n2-\delta\alpha-s}\left(\lambda_q^{s+\alpha}\|b_q\|_2\right)^{\delta}\left(\lambda_q^{s}\|b_q\|_2\right)^{3-\delta}\\
\leq &C\sum_{q\geq -1}\left(\lambda_q^{s+\alpha}\|b_q\|_2\right)^{\delta}\left(\lambda_q^{s}\|b_q\|_2\right)^{3-\delta}\\
\leq & \frac{\mu}{16}\sum_{q\geq -1}\lambda_q^{2s+2\alpha}\|b_q\|_2^2+C_\mu\sum_{q\geq -1}\left(\lambda_q^{s}\|b_q\|_2\right)^{\frac{2(3-\delta)}{2-\delta}}
\end{split}
\end{equation}
for $0<\delta<2$ and $s\geq 2+\frac n2-\delta\alpha>2+\frac n2-2\alpha$.

We continue to $I_{52}$ and decompose it by adapting commutator (\ref{comm-v2}),
\begin{equation}\notag
\begin{split}
I_{52}=&\sum_{q\geq-1}\sum_{|q-p|\leq 2}\lambda_q^{2s}\int_{\mathbb R^3}\Delta_q(\nabla\times b_{\leq p-2}\times b_{p})\cdot\nabla\times b_q\, dx\\
=&\sum_{q\geq-1}\sum_{|q-p|\leq 2}\lambda_q^{2s}\int_{\mathbb R^3}[\Delta_q, \nabla\times b_{\leq p-2}\times]b_p\cdot\nabla\times b_q\, dx\\
&+\sum_{q\geq-1}\lambda_q^{2s}\int_{\mathbb R^3}\nabla\times b_{\leq q-2} \times b_q\cdot\nabla\times b_q\, dx\\
&+\sum_{q\geq-1}\sum_{|p-q|\leq 2}\lambda_q^{2s}\int_{\mathbb R^3}\nabla\times (b_{\leq p-2}-b_{\leq q-2})\times (b_p)_q\cdot\nabla\times b_q\, dx\\
=&I_{521}+I_{522}+I_{523}.
\end{split}
\end{equation}
We will only show the estimate of $I_{522}$,  since $I_{521}$ enjoys the same estimate as $I_{511}$ due to the commutator estimate in Lemma \ref{le-Hall2} and $I_{523}$  can be estimated as $I_{513}$. 
Integration by parts,  identity (\ref{vec-ident1})
along with the fact that $\nabla\cdot b_q=0$ infers
\begin{equation}\notag
\begin{split}
I_{522}=&\sum_{q\geq-1}\lambda_q^{2s}\int_{\mathbb R^3}\nabla\times\left(\nabla\times b_{\leq q-2} \times b_q\right)\cdot b_q\, dx\\
=&\sum_{q\geq-1}\lambda_q^{2s}\int_{\mathbb R^3}\left[(b_q\cdot\nabla)\nabla\times b_{\leq q-2}-(\nabla\cdot \nabla\times b_{\leq q-2})b_q\right] \cdot b_q\, dx\\
&-\sum_{q\geq-1}\lambda_q^{2s}\int_{\mathbb R^3}(\nabla\times b_{\leq q-2}\cdot\nabla)b_q \cdot b_q\, dx.
\end{split}
\end{equation}
Since $\nabla\cdot (\nabla\times b_{\leq q-2})=0$, it is obvious the last integral vanishes.
Thus we have
\begin{equation}\notag
\begin{split}
|I_{522}|\leq &\sum_{q\geq-1}\lambda_q^{2s}\int_{\mathbb R^3}\left|\left[(b_q\cdot\nabla)\nabla\times b_{\leq q-2}-(\nabla\cdot \nabla\times b_{\leq q-2})b_q\right] \cdot b_q\right|\, dx\\
\lesssim &\sum_{q\geq-1}\lambda_q^{2s}\|\nabla^2 b_{\leq q-2}\|_\infty\|b_q\|_2^2\\
\lesssim &\sum_{q\geq-1}\lambda_q^{2s}\|b_q\|_2^2\sum_{p\leq q}\lambda_p^{2+\frac n2}\|b_p\|_2
\end{split}
\end{equation}
which share the same estimate of $I_{511}$. 

The last term $I_{53}$ is treated as
\begin{equation}\notag
\begin{split}
|I_{53}|\leq &\sum_{q\geq -1}\sum_{p\geq q-2}\lambda_q^{2s}\int_{\mathbb R^3}|\Delta_q(b_p\times \nabla\times\tilde b_p)\cdot\nabla\times b_q|\, dx\\
\lesssim &\sum_{q\geq -1}\lambda_q^{2s}\|\nabla b_q\|_\infty\sum_{p\geq q-3}\|b_p\|_2\|\nabla \tilde b_p\|_2\\
\lesssim &\sum_{q\geq -1}\lambda_q^{2s+1+\frac{n}2}\|b_q\|_2\sum_{p\geq q-3}\lambda_p\|b_p\|_2^2\\
\lesssim &\sum_{p\geq -1}\lambda_p\|b_p\|_2^2\sum_{q\leq p+3} \lambda_q^{2s+1+\frac{n}2}\|b_q\|_2
\end{split}
\end{equation}
which turns out to be similar as $I_{511}$ again. 
Summarizing the analysis above, we obtain
\begin{equation}\label{est-i5}
|I_5|\lesssim \frac\nu{8}\sum_{q\geq -1}\lambda_q^{2s+2\alpha}\|b_q\|_2^2+C_\mu\left(\sum_{q\geq -1}\lambda_q^{2s}\|b_q\|_2^2\right)^{\gamma_1}+C_\mu\left(\sum_{q\geq -1}\lambda_q^{2s}\|b_q\|_2^2\right)^{\gamma_2}
\end{equation}
for some $\gamma_1,\gamma_2>1$.
Putting together of (\ref{ineq-ubq}), (\ref{est-i3}), (\ref{est-i2}), and (\ref{est-i5}), there exist constants $C_\nu$, $C_\mu$, and $C_{\nu,\mu}$ such that
\begin{equation}\label{energy3}
\begin{split}
&\frac{d}{dt}\left(\| u\|_{\dot H^{s}}^2+\|b\|_{\dot H^{s}}^2\right)+\nu\sum_{q\geq -1}\lambda_q^{2s+2}\|u_q\|_2^2+\mu\sum_{q\geq -1}\lambda_q^{2s+2\alpha}\|b_q\|_2^2\\
\leq & C_\nu\left(\sum_{q\geq -1}\lambda_q^{2s}\|u_q\|_2^2\right)^{\gamma_1}+C_\nu\left(\sum_{q\geq -1}\lambda_q^{2s}\|u_q\|_2^2\right)^{\gamma_2}\\
&+C_\mu\left(\sum_{q\geq -1}\lambda_q^{2s}\|b_q\|_2^2\right)^{\gamma_1}+C_\mu\left(\sum_{q\geq -1}\lambda_q^{2s}\|b_q\|_2^2\right)^{\gamma_2}\\
\leq &C_{\nu,\mu}\left(\| u\|_{\dot H^{s}}^2+\|b\|_{\dot H^{s}}^2\right)^{\gamma_1}+C_{\nu,\mu}\left(\| u\|_{\dot H^{s}}^2+\|b\|_{\dot H^{s}}^2\right)^{\gamma_2}
\end{split}
\end{equation}
Notice that $\gamma_1, \gamma_2>1$ and hence the energy inequality (\ref{energy3}) is in the type of Riccati. It follows that, there exists a time $T>0$ which depends on $\nu,\mu$ and $\|u_0\|_{H^s}, \|b_0\|_{H^s}$ such that
\begin{equation}\notag
\|u(t)\|_{H^s}^2+\|b(t)\|_{H^s}^2\leq C(\nu,\mu,T, \|u_0\|_{H^s}, \|b_0\|_{H^s})\left(\|u_0\|_{H^s}^2+\|b_0\|_{H^s}^2\right)
\end{equation}
for $0\leq t<T$, and a constant $C$ depending on $\nu,\mu, T$ and $\|u_0\|_{H^s}, \|b_0\|_{H^s}$.

\cbdu



\section{Convergence of the Hall-MHD to the MHD system}
\label{sec:con}

In this section, we show that solutions $(u^\eta, b^\eta, p^\eta)$ of (\ref{HMHD}) with $\alpha=1$ in $H^{\frac n2}$ converges to a solution $(u, b, p)$ of the MHD system, as $\eta\to 0$.  Namely, we prove
\begin{Theorem}
Let $(u^\eta, b^\eta, p^\eta)$ be a solution to (\ref{HMHD}) with $\alpha=1$ obtained in Theorem \ref{thm} associated with initial data $(u_0,b_0)$. Let $(u, b, p)$ be a solution to (\ref{HMHD}) with $\eta=0$ and $\alpha=1$ under the same initial data. Then we have
\[\lim _{\eta\to 0}(\|u^\eta-u\|_2+\|b^\eta-b\|_2)=0.\]
\end{Theorem}
\pf
Take the difference $U=u^\eta-u$, $B=b^\eta-b$ and $\pi=p^\eta-p$, which satisfy the equations:
\begin{equation}\label{Eq-diff}
\begin{split}
U_t+u\cdot\nabla U-b\cdot\nabla B+U\cdot\nabla u^\eta-B\cdot\nabla b^\eta+\nabla \pi&=\nu\Delta U,\\
B_t+u\cdot\nabla B-b\cdot\nabla U+U\cdot\nabla b^\eta-B\cdot\nabla u^\eta-\eta\nabla\times ((\nabla\times b^\eta)\times b^\eta)&=\mu\Delta B,\\
\nabla \cdot U=0, \ \ \ \nabla\cdot B&=0.
\end{split}
\end{equation}
Multiplying the first equation by $U$ and the second by $B$,  we obtain (formally)
\begin{equation}\notag
\begin{split}
&\frac12\frac{d}{dt}\|U\|_2^2+\nu\|\nabla U\|_2^2\\
=&\int_{\R^3}b\cdot\nabla B\cdot U\, dx
-\int_{\R^3}U\cdot\nabla u^\eta\cdot U\, dx+\int_{\R^3}B\cdot\nabla b^\eta\cdot U\, dx,\\
&\frac12\frac{d}{dt}\|B\|_2^2+\mu\|\nabla B\|_2^2\\
=&\int_{\R^3}b\cdot\nabla U\cdot B\, dx
-\int_{\R^3}U\cdot\nabla b^\eta\cdot B\, dx+\int_{\R^3}B\cdot\nabla u^\eta\cdot B\, dx\\
&+\eta\int_{\R^3}\nabla\times((\nabla\times b^\eta)\times b^\eta)\cdot B\, dx.
\end{split}
\end{equation}
Adding the two yields, provided that $(u^\eta, b^\eta, p^\eta)$ and $(u, b, p)$ are regular enough,
\begin{equation}\notag
\begin{split}
&\frac12\frac{d}{dt}\left(\|U\|_2^2+\|B\|_2^2\right)+\nu\|\nabla U\|_2^2+\mu\|\nabla B\|_2^2\\
=&-\int_{\R^3}U\cdot\nabla u^\eta\cdot U\, dx+\int_{\R^3}B\cdot\nabla b^\eta\cdot U\, dx
-\int_{\R^3}U\cdot\nabla b^\eta\cdot B\, dx\\
&+\int_{\R^3}B\cdot\nabla u^\eta\cdot B\, dx
+\eta\int_{\R^3}\nabla\times((\nabla\times b^\eta)\times b^\eta)\cdot B\, dx\\
\equiv & I_1+I_2+I_3+I_4+I_5.
\end{split}
\end{equation}
It is straight forward to notice that 
\[|I_1+I_2+I_3+I_4|\leq C \left(\|\nabla u^\eta\|_\infty+\|\nabla b^\eta\|_\infty\right)\left(\|U\|_2^2+\|B\|_2^2\right);\]
and also
\[
\begin{split}
|I_1+I_2+I_3+I_4|\leq & C(\nu^{-1}+\mu^{-1}) \left(\|u^\eta\|_\infty+\| b^\eta\|_\infty\right)\left(\|U\|_2^2+\|B\|_2^2\right)\\
&+\frac14\nu \|\nabla U\|_2^2+\frac14\mu \|\nabla B\|_2^2.
\end{split}\]
We estimate $I_5$ as
\begin{equation}\notag
\begin{split}
|I_5|=&\left|\eta \int_{\R^3}((\nabla\times b^\eta)\times b^\eta)\cdot \nabla\times B\, dx\right|\\
\leq &C\eta \|\nabla b^\eta\|_\infty\|b^\eta\|_2\|\nabla B\|_2\\
\leq &C\eta^2\mu^{-1} \|\nabla b^\eta\|_\infty^2\|b^\eta\|_2^2+\frac1{4}\mu\|\nabla B\|_2^2
\end{split}
\end{equation}
or as
\begin{equation}\notag
\begin{split}
|I_5|=&\left|\eta \int_{\R^3}((\nabla\times b^\eta)\times b^\eta)\cdot \nabla\times B\, dx\right|\\
\leq &C\eta \| b^\eta\|_\infty\|\nabla b^\eta\|_2\|\nabla B\|_2\\
\leq &C\eta^2\mu^{-1} \|b^\eta\|_\infty^2\|\nabla b^\eta\|_2^2+\frac1{4}\mu\|\nabla B\|_2^2 
\end{split}
\end{equation}
Combining the above estimates leads to, for $s>\frac n2$
\begin{equation}\notag
\frac{d}{dt}\left(\|U\|_2^2+\|B\|_2^2\right)
\leq C \left(\|U\|_2^2+\|B\|_2^2\right)+C\eta^2\mu^{-1}\|\nabla b^\eta\|_2^2,
\end{equation}
from which Gr\"onwall's inequality implies that 
\[\|U(t)\|_2^2+\|B(t)\|_2^2\leq C\eta^2\mu^{-1}+(\|U(0)\|_2^2+\|B(0)\|_2^2+C\eta^2\mu^{-1})e^{Ct}.\]
Note that $U(0)=B(0)=0$. Thus 
\[\lim _{\eta\to 0}(\|U(t)\|_2^2+\|B(t)\|_2^2)=0,\]
and the convergence rate is $\mathcal O(\eta^2)$.
\cbdu

\bigskip

\section{Appendix}

\subsection{Littlewood-Paley decomposition}
\label{sec:LPD}
Our analysis is built on the Littlewood-Paley decomposition theory.
Basic languages and concepts are introduced briefly below.

We choose a nonnegative radial function $\chi\in C_0^\infty(\R^n)$ satisfying
\begin{equation}\notag
\chi(\xi)=
\begin{cases}
1, \ \ \mbox { for } |\xi|\leq\frac{3}{4}\\
0, \ \ \mbox { for } |\xi|\geq 1.
\end{cases}
\end{equation}
Denote $\lambda_q=2^q$ for integers $q$. A sequence of cut-off functions are defined,
\bg\notag
\varphi(\xi)=\chi(\frac{\xi}{2})-\chi(\xi), \ \ 
\varphi_q(\xi)=
\begin{cases}
\varphi(\lambda_q^{-1}\xi)  \ \ \ \mbox { for } q\geq 0,\\
\chi(\xi) \ \ \ \mbox { for } q=-1.
\end{cases}
\end{equation}
For a tempered distribution vector field $u$ we define the Littlewood-Paley projection
\begin{equation}\notag
\begin{split}
&h=\mathcal F^{-1}\varphi, \qquad \tilde h=\mathcal F^{-1}\chi,\\
&u_q:=\Delta_qu=\mathcal F^{-1}(\varphi(\lambda_q^{-1}\xi)\mathcal Fu)=\lambda_q^n\int h(\lambda_qy)u(x-y)dy,  \qquad \mbox { for }  q\geq 0,\\
& u_{-1}=\mathcal F^{-1}(\chi(\xi)\mathcal Fu)=\int \tilde h(y)u(x-y)dy,
\end{split}
\end{equation}
where $\mathcal F$ and $\mathcal F^{-1}$ denote the Fourier transform and inverse Fourier transform, respectively.
Due to the Littlewood-Paley theory, the identity
\bg\notag
u=\sum_{q=-1}^\infty u_q
\ed
holds in the sense of distribution, which is the fundamental idea of shell decomposition.
We also denote the various summation terms simply by
\bg\notag
u_{\leq Q}=\sum_{q=-1}^Qu_q, \qquad u_{(Q, N]}=\sum_{p=Q+1}^N u_p, \qquad \tilde u_q=\sum_{|p-q|\leq 1}u_p.
\ed

We can adapt the norm of Sobolev space $\dot H^s$ as
\[
  \|u\|_{\dot H^s} \sim \left(\sum_{q=-1}^\infty\lambda_q^{2s}\|u_q\|_2^2\right)^{1/2}, \ \ \ s\in\R.
\]

Bernstein's inequality satisfied by the dyadic blocks $u_q$ is introduced below.
\begin{Lemma}\label{le:bern} 
Let $n$ be the space dimension and $r\geq s\geq 1$. Then for all tempered distributions $u$, we have
\bg\notag
\|u_q\|_{r}\lesssim \lambda_q^{n(\frac{1}{s}-\frac{1}{r})}\|u_q\|_{s}.
\ed
\end{Lemma}

\bigskip


\subsection{Bony's paraproduct and commutators}
\label{sec-commu}

We adapt the following version of Bony's paraproduct 
\begin{equation}\label{Bony}
\begin{split}
\Delta_q(u\cdot\nabla v)=&\sum_{|q-p|\leq 2}\Delta_q(u_{\leq{p-2}}\cdot\nabla v_p)+
\sum_{|q-p|\leq 2}\Delta_q(u_{p}\cdot\nabla v_{\leq{p-2}})\\
&+\sum_{p\geq q-2} \Delta_q(\tilde u_p \cdot\nabla v_p),
\end{split}
\end{equation}
which is used through the paper to decompose the nonlinear terms. 
We introduce a commutator as
\begin{equation} \label{commu}
[\Delta_q, u_{\leq{p-2}}\cdot\nabla]v_p=\Delta_q(u_{\leq{p-2}}\cdot\nabla v_p)-u_{\leq{p-2}}\cdot\nabla \Delta_qv_p.
\end{equation}
\begin{Lemma}\label{le-commu}
The following estimate holds, for any $1<r<\infty$
\[\|[\Delta_q,u_{\leq{p-2}}\cdot\nabla] v_p\|_{r}\lesssim \|\nabla u_{\leq p-2}\|_\infty\|v_p\|_{r}.\]
\end{Lemma}

To treat the Hall term,  we recall a fundamental identity for vector valued functions $F$ and $G$,
\begin{equation}\label{vec-ident1}
\nabla \times(F\times G)=[(G\cdot\nabla) F-(\nabla\cdot F)G]-[(F\cdot\nabla) G-(\nabla\cdot G)F].
\end{equation}
In addition, two more commutators are defined
\begin{equation}\label{comm-v}
[\Delta_q,F\times\nabla\times]G=\Delta_q(F\times(\nabla\times G))-F\times(\nabla\times G_q),
\end{equation}
\begin{equation}\label{comm-v2}
[\Delta_q,(\nabla\times F)\times]G=\Delta_q((\nabla\times F)\times G)-(\nabla\times F)\times G_q.
\end{equation}
They satisfy the estimates below. 
\begin{Lemma}\label{le-Hall1}
Assume $\nabla\cdot F=0$ and $F$, $G$ vanish at large $|x|\in \R^3$. For any $1\leq r\leq \infty$, we have
\[\|[\Delta_q,F\times\nabla\times]G\|_r\lesssim  \|\nabla F\|_\infty\|G\|_r;\]
\[\|[\Delta_q, (\nabla\times F)\times]G\|_r\lesssim  \|\nabla F\|_\infty\|G\|_r.\]
\end{Lemma}

\begin{Lemma}\label{le-Hall2}
Assume the vector valued functions $F$, $G$ and $H$ vanish at large $|x|\in \R^3$. For any $1\leq r_1, r_2\leq \infty$ with $\frac1{r_1}+\frac1{r_2}=1$,  we have
\[
\left|\int_{\R^3}[\Delta_q, (\nabla\times F)\times]G\cdot\nabla\times H\, dx\right|
\lesssim \|\nabla^2 F\|_\infty\|G\|_{r_1}\|H\|_{r_2}.
\]
\end{Lemma}


\begin{thebibliography}{XX}



\bibitem{ADFL}
M. Acheritogaray, P. Degond, A. Frouvelle and J-G. Liu.
\newblock {\em Kinetic formulation and global existence for the Hall-Magnetohydrodynamic system}.
\newblock Kinetic and Related Models, 4: 901--918, 2011.


\bibitem{CDL}
D. Chae, P. Degond and J-G. Liu.
\newblock {\em Well-posedness for Hall--magnetohydrodynamics}.
\newblock arXiv:1212.3919, 2012.

\bibitem{CL}
D. Chae and J. Lee.
\newblock {\em On the blow-up criterion and small data global existence for the Hall-magneto-hydrodynamics}.
\newblock J. Differential Equations, 256: 3835--3858, 2014.

\bibitem{CS}
D. Chae,  and M. Schonbek.
\newblock {\em On the temporal decay for the Hall-magnetohydrodynamic equations}.
\newblock arXiv:1302.4601, 2013.

\bibitem{CWW}
D. Chae,  R. Wan and J. Wu.
\newblock {\em Local well-posedness for the Hall--MHD equations with fractional magnetic diffusion}.
\newblock arXiv:1404.0486v2, 2014.

\bibitem{CWeng}
D. Chae and S. weng.
\newblock {\em Singularity formation for the incompressible Hall-MHD equations without resistivity}.
\newblock Ann. I. H. Poincar\'e-AN, Vol. 33: 1009--1022, 2016.

\bibitem{CW}
D. Chae and J. Wolf.
\newblock {\em On partial regularity for the 3D non-stationary Hall magnetohydrodynamics equations on the plane}.
\newblock arXiv:1502.0347, 2015.






\bibitem{D}
M. Dai.
\newblock {\em Regularity criterion and energy conservation for the supercritical Quasi-Geostrophic equation}.
\newblock Journal of Mathematical Fluid Mechanics. To appear. ArXiv:1505.02293, 2015.

\bibitem{DL}
M. Dai and H. Liu.
\newblock {\em Long time behavior of solutions to the 3D Hall-magneto-hydrodynamics system with one diffusion}.
\newblock arXiv:1705.02647, 2017.

\bibitem{DS}
E. Dumas and F. Sueur.
\newblock {\em On the weak solutions to the Maxwell-Landau-Lifshitz equations and to the Hall-magnetohydrodynamic equations}.
\newblock Comm. Math. Phys., 330: 1179--1225, 2014.

\bibitem{FMRR}
C. L. Fefferman, D. S. McCormick, J. C. Robinson and J. L. Rodrigo.
\newblock {\em Higher order commutator estimates and local existence for the non-resistive MHD equations and related models}.
\newblock Journal of Functional Analysis, Vol. 267: 1035--1056, 2014.



\bibitem{Gr}
L. Grafakos.
\newblock {\em Modern Fourier analysis}.
\newblock Second edition. Graduate Texts in Mathematics, 250. Springer, New York, 2009.













\end{thebibliography}
\end{document}